# Rényi and Shannon Entropies of Finite Mixtures of Multivariate Skew t-distributions

Salah H. Abid [1,*] and Uday J. Quaez [2]

[1] *Department of Mathematics, Education College, Al-Mustansiriya University, Baghdad, Iraq.* E-mail: [*] *abidsalah@uomustansiriyah.edu.iq.*

[2] *Department of Mathematics, Education College, Al-Mustansiriya University, Baghdad, Iraq.*

Shannon and Rényi entropies are quantitative measures of uncertainty in a data set. They are developed by Rényi in the context of entropy theory. These measures have been studied in the case of the multivariate t-distributions. We extend these tools to the class of multivariate skew t-distributions and then to more families of finite mixture of multivariate skew t[distributions. In particular, by using generalized HÖlder's inequality and some properties of multinomial theorem, we find upper and lower bounds of entropies for these families. An approximate value of these entropies can be calculated. In addition, an asymptotic expression for Rényi entropy is given by approximation and by using some inequalities and properties of $L^p$ -spaces. Finally, we give a real data examples to illustrate the behavior of entropy of the mixture model under consideration.

## 1. Introduction

Due to its flexibility, the finite mixture models have become widely used in the modeling and analysis of heterogeneous data sets. These models are an important statistical tool for many applications such as density estimation, data mining, medicine, image processing, satellite imaging and pattern recognition etc. (see for more detail [8], [9], [10], [11]). In multivariate analysis, [12] Azzalini A. and Dalla Valle A., introduced the multivariate skew normal distribution as an alternative to multivariate normal distribution to deal with skew in the data. In [16], Lin T., Lee J. and Wan H. proposed a development of mixture of skew normal and skew t models. Finite mixture of multivariate skew normal and skew t distributions were studied in Pyne S. [17]. More recently, Lee S. and McLachlan G. in [9] provided an overview of developments of finite mixture of skew t- distributions.

In other hand, Shanon in [18] proposed a measure to quantify the uncertainty of an event. In [19] Rényi generalized this measure for probability distribution which means the sensitive to the fine details of a density function. Javier E. and Contreras-Reyes J. [15] discussed the values of Rényi entropy of flexible class of skew normal distributions. Wood, Blythe and Evans [20] introduced some results for Rényi entropy of the totally asymmetric exclusion process. Also they calculated explicitly Rényi entropy whereby the squares of configuration probabilities are summed. It is important to note that the authors in [10] discussed the bounds and approximation of Rényi entropy of a class of mixture models of multivariate skew Gaussian by using multinomial theorem and generalized HÖlder's inequality.

In this paper, we propose a model of finite mixture of multivariate skew t-distributions. The explicitly expression of Shannon and Rényi entropies of skew t distribution is derived. By using generalized HÖlder's inequality and some properties of multinomial theorem, we find upper and lower bounds for Shannon and Rényi entropies of mixture models. An approximate value of these entropies can be calculated. In addition, an asymptotic expression for Rényi entropy is given by approximation and by using some inequalities and properties of $L^p$ -spaces. Finally, we give a real data examples to illustrate the behavior of Rényi entropy with the parameters $\alpha$ , skewness and freedom degrees of the proposed mixture model.

The remainder of this paper is organized as follow: In section 2. we begin with a preliminary material of mixture models and the measures (Shannon and Rényi) of information. Section 3. provides a description of asymptotic expression for Shannon and Rényi entropies of multivariate skew t- distributions by using some methods of numerical integration such as Monte Carlo and importance sampling methods. By using some theorems which related with multinomial theorem and generalized HÖlder's inequality, in section 4. we find upper and lower bounds and also we study an approximate Shannon and Rényi entropies for a finite mixture of multivariate skew t- distributions.

## 2. Preliminary Material

In this section, we introduce some basic definitions, notations and lemmas related to entropy theory and mixture models of skew distributions that we shall need them later on.

multivariate skew t-distribution has been proposed by Azzalini and Capitaino [7]. A random vector $X \in R^d$ has a multivariate skew t-distribution with location vector $\mu \in R^d$, scale matrix $S \in R^{d \times d}$, skewness vector $\delta \in R^d$ and freedom degrees $v \in R$, denoted by $(X \sim MST_d(\mu, S, \delta, v))$, if its density is given by

$$f_d(x; \mu, S, \delta, v) = 2g_d(x; \mu, S, v). G_1\left(\frac{1}{\beta}\hat{\delta}'S^{-1}(x-\mu)\sqrt{\frac{v+d}{v+(x-\mu)'S^{-1}(x-\mu)}}; 0,1, v+d\right) \tag{1}$$

where,

$g_d(:; \mu, S, v)$ is a d-dimensional multivariate t-distribution and $G_1$ is the distribution function of univariate standard t-distribution with freedom degrees v+d.

$\beta = \sqrt{1 - \hat{\delta}' S^{-1} \hat{\delta}} < 1$

$\hat{\delta} = \frac{1}{(1+\delta'\tilde{S}^{-1}S\tilde{S}^{-1}\delta)^{\frac{1}{2}}} S\tilde{S}^{-1}\delta$

$\tilde{S} = \text{diag}(s_{11}, s_{22}, \ldots, s_{dd})^{\frac{1}{2}}$, when $S=(s_{ij})$, $i, j = 1,2,\ldots,d$

The stochastic representation of X can be written as

$X = \mu + \hat{\delta}|Z_1| + Z_d$

where, $Z_1 \sim N_1\left(0, \frac{1}{w}\right)$, $Z_d \sim MN_d\left(0, \frac{1}{w}S\right)$ which are independent, the notations $N_1$, $MN_d$ and $|Z_1|$ are represent univariate, multivariate normal distributions and the absolute value of $Z_1$ respectively. It can be shown that the density function of X under this representation is given by equation (1). Lee and McLachlan [9] show that the multivariate skew t-distribution is multivariate skew normal distribution as $v \to \infty$ and multivariate t-distribution when $\delta \to 0$. The mean vector and covariance matrix of X are derived by Azzalini and Capitaino [7] in the following form

$$E(X) = \mu + \frac{\Gamma\left(\frac{v-1}{2}\right)\sqrt{\frac{v}{\pi}}}{\Gamma\left(\frac{v}{2}\right)} \hat{\delta} \tag{2}$$

$$Var(X) = \frac{v}{v-2}S - \left(\frac{\Gamma\left(\frac{v-1}{2}\right)}{\Gamma\left(\frac{v}{2}\right)}\right)^2 \frac{v}{\pi}\hat{\delta}\hat{\delta}' \tag{3}$$

**Definition 1.** Let X be a d-dimensional random vector which comes from an m-component mixtures of multivariate skew t-distributions. Then the density function of $X \sim FMMST_d(\mu, S, \delta, v, \varepsilon)$ is given by the following form

$$f(x; \mu, S, \delta, v, \varepsilon) = \sum_{i=1}^{n} \varepsilon_i f(x; \mu_i, S_i, \delta_i, v_i) \tag{4}$$

where, $\varepsilon_i$ denotes the mixing probability with $\varepsilon_i \geq 0$, $\sum_{i=1}^{n} \varepsilon_i = 1$, $f(x; \mu_i, S_i, \delta_i, v_i)$ represent the pdf of an m-component mixture model with parameter vector set $(\mu, S, \delta)$; $\mu = \{\mu_1, \mu_2, \ldots, \mu_m\}$ a set of vectors represent location parameters, $S = \{S_1, S_2, \ldots, S_m\}$ a set of dispersion matrices, the shape vector parameter is $\delta = \{\delta_1, \delta_2, \ldots, \delta_m\}$ and $v = \{v_1, v_2, \ldots, v_m\}$ is a set of freedom degrees.

If $\kappa = (\kappa_1, \kappa_2, \ldots, \kappa_n)$ represent a set of n latent allocations for densities of observations x then $f(x; \mu, S, \delta, v, \varepsilon) = \prod_{j=1}^{n} f(x; \mu, S, \delta, v, \kappa_j)$, where $\Pr(\kappa_j = i | \varepsilon) = \varepsilon_i$ then for any j-th component density in (4) is obtained as

$X_j|(\kappa_j = i) \stackrel{d}{=} \mu_i + \hat{\delta}_i|Z_{1j}| + Z_{dj}$, $j = 1, 2, \ldots, n$, where $Z_{1j} \sim N_1\left(0, \frac{1}{w_j}\right)$, $Z_d \sim N_d\left(0, \frac{1}{w_j}S_j\right)$ which are mutually independent and $\hat{\delta}_i = \frac{1}{(1+\delta_i'\tilde{S}_i^{-1}S_i\tilde{S}_i^{-1}\delta_i)^{\frac{1}{2}}} S_i\tilde{S}_i^{-1}\delta_i$, $i = 1, 2, \ldots, m$. Equations (2-3) gives the first and second moments for each i-th component of X respectively, the mean and variance of X can be obtained as follow

$$E(X) = \sum_{i=1}^{m} \varepsilon_i \left(\mu_i + \frac{\Gamma\left(\frac{v_i-1}{2}\right)\sqrt{\frac{v_i}{\pi}}}{\Gamma\left(\frac{v_i}{2}\right)} \hat{\delta}_i\right) \tag{5}$$

$$Var(X) = \sum_{i=1}^{m} \varepsilon_i \frac{v_i}{v_i-2} S_i - \left(\sum_{i=1}^{m} \varepsilon_i \left(\frac{\Gamma\left(\frac{v_i-1}{2}\right)\sqrt{\frac{v_i}{\pi}}}{\Gamma\left(\frac{v_i}{2}\right)} \hat{\delta}_i\right)\right)\left(\sum_{i=1}^{m} \varepsilon_i \left(\frac{\Gamma\left(\frac{v_i-1}{2}\right)\sqrt{\frac{v_i}{\pi}}}{\Gamma\left(\frac{v_i}{2}\right)} \hat{\delta}_i\right)\right)' \tag{6}$$

**Definition 2.** Let X be a continuous random vector in $R^d$ with probability density function $f(x; \theta)$. Then the Shannon differential entropy is defined as

$$H(X;\theta) = -E(\ln(f(x;\theta))) \tag{7}$$

**Definition 3.** An αth-order Rényi differential entropy of a continuous random vector $X \in R^d$ with probability density function $f(x;\theta)$ is defined as

$$R_\alpha(X;\theta) = \begin{cases} \frac{1}{1-\alpha} \ln(E(f(x;\theta))^{\alpha-1}), & 0 < \alpha < \infty, \alpha \neq 1 \\ -E(\ln(f(x;\theta))) & \alpha = 1 \end{cases} \tag{8}$$

The relationship between Shannon and Rényi entropies is obtained by the limit $H(X;\theta) = \lim_{\alpha \to 1} R_\alpha(X;\theta)$. Translation does not change the entropy $H(X + C) = H(X) + C$ where, C is constant. Also for any $0 < \alpha_1 < \alpha_2$ the important property in Rényi entropy gives $R_{\alpha_2}(X;\theta) < R_{\alpha_1}(X;\theta)$ (see, e.g.,[8]).

**Definition 4.** [14] The digamma function is defined as the natural logarithmic derivative of gamma function

$$\psi(v) = \frac{d}{dv} \ln(\Gamma(v)) = \frac{\frac{d}{dv}(\Gamma(v))}{\Gamma(v)} \tag{9}$$

**Lemma 1.** [12] If a random vector $X \in R^d$ has zero mean and covariance matrix $R = E(XX')$, then the following inequality is accomplished
$$H(X;\theta) \leq \frac{1}{2} \ln(\det(2\pi \exp(1) R)) \tag{10}$$
with equality if and only if $X \sim MN_d(0, S)$.

**Lemma 2.** [11] Consider $y \sim MST_d(0, I, \delta, v)$ which defined in equation (1). Then

$$\delta' y \sqrt{\frac{v+d}{v+y'y}} \stackrel{d}{=} \frac{\sqrt{(v+d)\delta'\delta}\, x}{\sqrt{v+d-1+x^2}} \tag{11}$$

where, $x \sim ST_1(0, 1, \sqrt{\delta'\delta}, v + d - 1)$ and $\stackrel{d}{=}$ denotes equality in distribution.

**Lemma 3.** [5] Let $g_d$ be a d-dimensional probability density function centrally symmetric about 0, $G_1(.)$ be a continuous distribution function on the real line, so that $g_d(x) = \frac{dG_1(x)}{dx}$ exists a.e. and is an even density function and an odd real-valued function on $R^d$ by $\vartheta$. Thus, $f_d(x) = 2g_d(x). G_1(\vartheta(x))$ is a density function.

**Proposition 4.** [5] If the random vectors $X_0, X \in R^d$ have densities $g_d$ and $f_d$ respectively, where $g_d$ and $f_d$ satisfies the conditions in lemma 3., then $h(x_0) \stackrel{d}{=} h(x)$ for any even d-dimensional function h on $R^d$, irrespectively of factor $F(x) = G_1(\vartheta(x))$ where, F is distribution function.

**Lemma 5.** [13] **(Multinomial Theorem)** Let $\varepsilon_1, \varepsilon_2, \ldots, \varepsilon_n \geq 0$ and $r_1, r_2, \ldots, r_n \geq 0$ be any non-negative real numbers. Then the equality

$$\left(\sum_{i=1}^n \varepsilon_i r_i\right)^\alpha = \sum_{k_i \in B} \frac{\alpha!}{\prod_{i=1}^n k_i!} \prod_{i=1}^n (\varepsilon_i r_i)^{k_i} \tag{12}$$

is satisfied under the condition $\sum_{k_i \in B} \frac{\alpha!}{\prod_{j=1}^n k_j!} = n^\alpha$ where, $B = \{k_i \in N, \sum_{i=1}^n k_i = \alpha, i = 1, 2, \ldots, n\}$

**Lemma 6.** [13] Let $a_1, a_2, \ldots, a_n$ and $x_1, x_2, \ldots, x_n$ be two arbitrary sets of real numbers. If $\alpha$ is a positive integer then the following equality

$$\left(\sum_{i=1}^n a_i x_i\right)^\alpha = \left(\sum_{i=1}^n a_i\right)^\alpha x_n^\alpha + \sum_{i=1}^{n-1} \left(\sum_{k=1}^i a_k\right)^\alpha (x_i^\alpha - x_{i+1}^\alpha) + \sum_{k_t \in A} \frac{\alpha!}{\prod_{j=1}^n k_j!} \left(\prod_{t=1}^i (\varepsilon_t)^{k_t}\right) \left[\left(\prod_{t=1}^{i-1} x_t^{k_t}\right) - x_i^{\alpha - k_i}\right] \tag{13}$$

is accomplished. The set A is defined as $A = \{k_t \in N; 0 < k_t < \alpha, \sum_{i=1}^n k_i = \alpha, k_{t+1} = k_{t+2} = \cdots = k_n = 0\}$

**Lemma 7.** [13] Let $\varepsilon_1, \varepsilon_2, \ldots, \varepsilon_n \geq 0$ and $r_1, r_2, \ldots, r_n \geq 0$ be given. Then for any real numbers $p \geq 0$ and $0 \leq \alpha \leq p$, the following inequality is holds:

$$\left(\sum_{i=1}^n \varepsilon_i r_i\right)^\alpha \|r\|_p^{p-\alpha} \geq \sum_{i=1}^{n-1} (i)^{1-\frac{\alpha}{p}} \left(\sum_{k=1}^i \varepsilon_k\right)^\alpha (r_i^p - r_{i+1}^p) + n^{1-\frac{\alpha}{p}} \left(\sum_{k=1}^n \varepsilon_k\right)^\alpha r_n^p \tag{14}$$

where, $\|r\|_p = (\sum_{k=1}^{n} r_k^p)^{\frac{1}{p}}$

**Proposition 8.**
Let $X_0 \sim MT_d(\mu, S, v)$. Then the Shannon and Rényi entropies of $X_0$ are given in the following forms respectively

$$H(x_0; \mu, S, v) = \ln\left(\frac{\Gamma(\frac{v}{2})(v\pi)^{\frac{d}{2}}}{\Gamma(\frac{v+d}{2})}\right) + \frac{1}{2}\ln(\det(S)) + \left(\frac{v+d}{2}\right)\{\psi(v+d) - \psi(v)\} \tag{15}$$

$$R_\alpha(x_0; \mu, S, v) = H(x_0; \mu, S, v) + \frac{1}{(1-\alpha)}\ln\left(\frac{\Gamma(\frac{v+d}{2})\Gamma(\frac{\alpha v + (\alpha-1)d}{2})}{\Gamma(\frac{v}{2})\Gamma(\frac{\alpha(v+d)}{2})}\right), \; 0 < \alpha < \infty, \; \alpha \neq 1 \tag{16}$$

## 3. Shannon and Rényi Entropies for Multivariate Skew t- Distributions

In this section, we will derived explicit expressions for Shannon and Rényi entropies for multivariate skew t-distributions by using some properties of transformation and integrations. Also we give a simple illustrative example explains the relationship between the parameters α , δ and v with Shannon and Renyi entropies by using methods of numerical integration such as Monte Carlo and importance sampling methods.

**Proposition 9.** Let $X \sim MST_d(\mu, S, \delta, v)$ and let $Y \in R^d$ has a multivariate skew normal distribution with location vector $\mu \in R^d$, scale matrix $S \in R^{d \times d}$ and skewness vector $\delta \in R^d$ denoted by $Y \sim MSN_d(\mu, S, \delta)$. Then
  i. $H(Y; \mu, S, \delta) = \lim_{v \to \infty} H(X; \mu, S, \delta, v)$
  ii. $R_\alpha(Y; \mu, S, \delta) = \lim_{v \to \infty} R_\alpha(X; \mu, S, \delta, v)$

**Proof**
Let $f_v(X; \mu, S, \delta) = f(X; \mu, S, \delta, v)$, $v = 1, 2, \ldots$ be a sequence of multivariate skew t-distributions. Then
$\lim_{v \to \infty} f_v(X; \mu, S, \delta) = f(Y; \mu, S, \delta)$, where $Y \sim MSN_d(\mu, S, \delta)$
Applying Lebesgue dominated convergence theorem, to obtain that
$$\lim_{v \to \infty} E(\ln(f_v(X; \mu, S, \delta)))) = E(\ln(f(Y; \mu, S, \delta))))$$
In the similar method to prove part ii. This proposition shows that our problem formulation is generalized by the work of Contreras-Reyes and Cortés [10].

**Lemma 10.** Consider $X \sim MST_d(\mu, S, \delta, v))$. Then

$$E\left\{\ln\left(2G_1\left(\frac{1}{\beta}\hat{\delta}'S^{-1}(x-\mu)\sqrt{\frac{v+d}{v+(x-\mu)'S^{-1}(x-\mu)}}; 0, 1, v+d\right)\right)\right\}$$
$$= E\left\{\ln\left(2G_1\left(\frac{\sqrt{(v+d)\tilde{\delta}'\tilde{\delta}}\, y}{\beta\sqrt{v+d-1+y^2}}; 0, 1, v+d\right)\right) * 2G_1\left(\frac{1}{\beta^2}\sqrt{\tilde{\delta}'\tilde{\delta}}\, y\sqrt{\frac{(v+1)}{v+y^2}}; 0, 1, v+d\right)\right\} \tag{17}$$

where, $y \sim T_1(0, 1, v + d - 1)$ and $\tilde{\delta} = S^{-\frac{1}{2}}\hat{\delta}$

**Proof:** we have directly that

$$E\left\{\ln\left(2G_1\left(\frac{1}{\beta}\hat{\delta}'S^{-1}(x-\mu)\sqrt{\frac{v+d}{v+(x-\mu)'S^{-1}(x-\mu)}}; 0, 1, v+d\right)\right)\right\}$$
$$= \int_{R^d} \ln\left(2G_1\left(\frac{1}{\beta}\hat{\delta}'S^{-1}(x-\mu)\sqrt{\frac{v+d}{v+(x-\mu)'S^{-1}(x-\mu)}}; 0, 1, v+d\right)\right) g_d(x; \mu, S, v)$$
$$* 2G_1\left(\frac{1}{\beta}\hat{\delta}'S^{-1}(x-\mu)\sqrt{\frac{v+d}{v+(x-\mu)'S^{-1}(x-\mu)}}; 0, 1, v+d\right) dx \tag{18}$$

Using the change of variables $z = S^{-\frac{1}{2}}(x - \mu)$ and $\tilde{\delta} = S^{-\frac{1}{2}}\hat{\delta}$ associated with Jacobian matrix is given by $\det(S)^{\frac{1}{2}}$, to get that $Z \sim MST_d(0, I, \delta, v))$. Hence

$$E\left\{\ln\left(2G_1\left(\frac{1}{\beta}\hat{\delta}'S^{-1}(x-\mu)\sqrt{\frac{v+d}{v+(x-\mu)'S^{-1}(x-\mu)}}; 0, 1, v+d\right)\right)\right\}$$
$$= \int_{R^d} \ln\left(2G_1\left(\frac{1}{\beta}\tilde{\delta}'z\sqrt{\frac{v+d}{v+z'z}}; 0, 1, v+d\right)\right) g_d(z; 0, I, v)$$

$$* 2G_1\left(\frac{1}{\beta}\tilde{\delta}'z\sqrt{\frac{v+d}{v+z'z}}; 0,1, v+d\right)dz \qquad (19)$$

Therefore,

$$E\left\{\ln\left(2G_1\left(\frac{1}{\beta}\hat{\delta}'S^{-1}(x-\mu)\sqrt{\frac{v+d}{v+(x-\mu)'S^{-1}(x-\mu)}}; 0,1, v+d\right)\right)\right\}$$

$$= E\left\{\ln\left(2G_1\left(\frac{1}{\beta}\tilde{\delta}'z\sqrt{\frac{v+d}{v+z'z}}; 0,1, v+d\right)\right)\right\} \qquad (20)$$

Lemma 2. implies that

$$E\left\{\ln\left(2G_1\left(\frac{1}{\beta}\hat{\delta}'S^{-1}(x-\mu)\sqrt{\frac{v+d}{v+(x-\mu)'S^{-1}(x-\mu)}}; 0,1, v+d\right)\right)\right\}$$

$$= E\left\{\ln\left(2G_1\left(\frac{1}{\beta}\frac{\sqrt{(v+d)\tilde{\delta}'\tilde{\delta}}\,w}{\sqrt{v+d-1+w^2}}; 0,1, v+d\right)\right)\right\} \qquad (21)$$

where, $w \sim ST_1(0,1, \sqrt{\tilde{\delta}'\tilde{\delta}}, v+d-1)$

**Proposition 11.** Let $X_0 \sim MT_d(\mu, S, v))$ and $X \sim MST_d(\mu, S, \delta, v)$. Then the Shannon entropy can be written as follows:

$$H(X; \mu, S, \delta, v) = H(x_0; \mu, S, v) - \mathfrak{S}_{v,d}(\tilde{\delta}) \qquad (22)$$

where,

$$\mathfrak{S}_{v,d}(\tilde{\delta}) = E\left\{\ln\left(2G_1\left(\frac{\sqrt{(v+d)\tilde{\delta}'\tilde{\delta}}\,y}{\beta\sqrt{v+d-1+y^2}}; 0,1, v+d\right)\right) * 2G_1\left(\frac{1}{\beta^2}\sqrt{\tilde{\delta}'\tilde{\delta}}\,y\sqrt{\frac{(v+1)}{v+y^2}}; 0,1, v+d\right)\right\} \qquad (23)$$

$y \sim T_1(0,1, v+d-1)$

**Proof:** by taking the natural logarithm and expectation for both sides of equation (1), we have
$E\left(\ln(f_d(x; \mu, S, \delta, v))\right)$

$$= E\left(\ln(g_d(x; \mu, S, v))\right) + E\left\{\ln\left(2G_1\left(\frac{1}{\beta}\hat{\delta}'S^{-1}(x-\mu)\sqrt{\frac{v+d}{v+(x-\mu)'S^{-1}(x-\mu)}}; 0,1, v+d\right)\right)\right\}$$

Using lemma 3. and proposition 4. , we obtain

$$H(X; \mu, S, \delta, v) = H(x_0; \mu, S, v) - E\left\{\ln\left(2G_1\left(\frac{1}{\beta}\hat{\delta}'S^{-1}(x-\mu)\sqrt{\frac{v+d}{v+(x-\mu)'S^{-1}(x-\mu)}}; 0,1, v+d\right)\right)\right\}$$

Lemma 10. gives us the required result of this proof.

**Lemma 12.** Suppose that $X \sim MST_d(\mu, S, \delta, v)$. Then:

$$\int_{R^d}(f(x; \mu, S, \delta, v))^\alpha dx = \hat{C}_{\alpha,v,d}(S) * E\left\{2\left(G_1\left(\frac{\sqrt{(v+d)\tilde{\delta}'\tilde{\delta}}\,x}{\beta\sqrt{\alpha(v+d)-1+x^2}}; 0,1, v+d\right)\right)^\alpha\right\} \qquad (24)$$

where,

$$\hat{C}_{\alpha,v,d}(S) = \left(\frac{\Gamma\left(\frac{v}{2}\right)}{\Gamma\left(\frac{v+d}{2}\right)}\right)^{1-\alpha}(v\pi)^{\frac{d}{2}(1-\alpha)}(\det(S))^{\frac{1}{2}(1-\alpha)}\frac{\Gamma\left(\frac{v+d}{2}\right)\Gamma\left(\frac{\alpha(v+d)-d}{2}\right)}{\Gamma\left(\frac{v}{2}\right)\Gamma\left(\frac{\alpha(v+d)}{2}\right)} \qquad (25)$$

$x \sim T_1(0,1, \alpha(v+d)-d), \qquad 0 < \alpha < \infty, \alpha \neq 1$

**Proof:** Replacing $\tilde{\delta}$ by $S^{-\frac{1}{2}}\hat{\delta}$ and transforming the variable $z = S^{-\frac{1}{2}}(x-\mu)$ associated with Jacobian matrix $S^{\frac{1}{2}}$ in equation (1), we get

$\int_{R^d}(f(x; \mu, S, \delta))^\alpha dx =$

$$\left(\frac{2\Gamma\left(\frac{v+d}{2}\right)(v\pi)^{\frac{-d}{2}}}{\Gamma\left(\frac{v}{2}\right)}\right)^\alpha (\det(S))^{\frac{1}{2}(1-\alpha)}\int_{R^d}\left(1+\frac{1}{v}(z'z)\right)^{-\alpha\left(\frac{v+d}{2}\right)}\left(G_1\left(\frac{1}{\beta}\tilde{\delta}'z\sqrt{\frac{v+d}{v+z'z}}; 0,1, v+d\right)\right)^\alpha dz$$

Again replacing $\alpha(v+d)-d$ by $u$ and transforming $\sqrt{\frac{u}{v}}z$ by the variable $y$, we obtain

$\int_{R^d}(f(x;\mu,S,\delta))^\alpha dx =$
$\left(\frac{2\Gamma\left(\frac{v+d}{2}\right)(v\pi)^{-\frac{d}{2}}}{\Gamma\left(\frac{v}{2}\right)}\right)^\alpha (\det(S))^{\frac{1}{2}(1-\alpha)}\left(\frac{v}{u}\right)^{\frac{d}{2}}\int_{R^d}\left(1+\frac{y'y}{u}\right)^{-\left(\frac{u+d}{2}\right)}\left(G_1\left(\frac{1}{\beta}\tilde{\delta}'y\sqrt{\frac{v+d}{u+y'y}};0,1,v+d\right)\right)^\alpha dy$

Consequently,

$\int_{R^d}(f(x;\mu,S,\delta))^\alpha dx =$
$=\left(\frac{2\Gamma\left(\frac{v+d}{2}\right)(v\pi)^{-\frac{d}{2}}}{\Gamma\left(\frac{v}{2}\right)}\right)^\alpha (\det(S))^{\frac{1}{2}(1-\alpha)}\left(\frac{v}{u}\right)^{\frac{d}{2}}\frac{\Gamma\left(\frac{u}{2}\right)(u\pi)^{\frac{d}{2}}}{\Gamma\left(\frac{u+d}{2}\right)}E\left\{\left(G_1\left(\frac{1}{\beta}\tilde{\delta}'y\sqrt{\frac{v+d}{u+y'y}};0,1,v+d\right)\right)^\alpha\right\}$

where, $y \sim T_d(0, I, u)$. This proof is completed by using lemma 2.

In the right side of equation(25), Computing the expectation in a convenient and fast method needs to be based on numerical integration.

**Corollary 13.** If $X \sim MST_d(\mu, S, \delta, v)$ and $X_0 \sim MT_d(\mu, S, v)$, then the Rényi and Shannon entropies may be written as

i. $R_\alpha(X;\mu,S,\delta,v) = R_\alpha(X_0;\mu,S,v) + \frac{1}{(1-\alpha)}\ln\left(E\left\{\left(2G_1\left(\frac{\sqrt{(v+d)}\tilde{\delta}'\tilde{\delta}\,x}{\beta\sqrt{\alpha(v+d)-1+x^2}};0,1,v+d\right)\right)^\alpha\right\}\right)$ (26)

ii. $\lim_{\alpha \to 1} R_\alpha(X;\mu,S,\delta,v) = H(X;\mu,S,\delta,v)$ (27)

**Proof :** Taking natural logarithm and multiplying by $\frac{1}{1-\alpha}$ for both sides of equation (24), we get
$R_\alpha(X;\mu,S,\delta,v)$

$= \ln\left(\frac{\Gamma\left(\frac{v}{2}\right)}{\Gamma\left(\frac{v+d}{2}\right)}\right) + \frac{d}{2}\ln(v\pi) + \frac{1}{2}\ln(\det(S)) + \frac{1}{(1-\alpha)}\ln\left(\frac{\Gamma\left(\frac{v+d}{2}\right)\Gamma\left(\frac{\alpha v+(\alpha-1)d}{2}\right)}{\Gamma\left(\frac{v}{2}\right)\Gamma\left(\frac{\alpha(v+d)}{2}\right)}\right)$
$+ \frac{1}{(1-\alpha)}\ln\left(E\left\{2\left(G_1\left(\frac{1}{\beta}\frac{\sqrt{(v+d)}\delta'\delta\,x}{\sqrt{u+d-1+x^2}};0,1,v+d\right)\right)^\alpha\right\}\right)$ (28)

where, $x \sim T_1(0,1,u)$
Equation (16) implies that

$R_\alpha(X;\mu,S,\delta,v) = R_\alpha(X_0;\mu,S,v) + \frac{1}{(1-\alpha)}\ln\left(E\left\{\left(2G_1\left(\frac{\sqrt{(v+d)}\tilde{\delta}'\tilde{\delta}\,x}{\beta\sqrt{\alpha(v+d)-1+x^2}};0,1,v+d\right)\right)^\alpha\right\}\right)$

$x \sim T_1(0,1,\alpha(v+d)-d)$
Now, to prove part ii. By taking the limiting as α converges to 1, we have
$\lim_{\alpha \to 1} R_\alpha(X;\mu,S,\delta,v)$

$= \lim_{\alpha \to 1} R_\alpha(X_0;\mu,S,v) + \lim_{\alpha \to 1}\frac{1}{(1-\alpha)}\ln\left(E\left\{\left(2G_1\left(\frac{1}{\beta};\frac{\sqrt{v+d}\delta'\delta\,x}{\sqrt{\alpha(v+d)-1+x^2}};0,1,v+d\right)\right)^\alpha\right\}\right)$

$= H(x_0;\mu,S,v) - E\left\{\ln\left(2G_1\left(\frac{\sqrt{(v+d)}\tilde{\delta}'\tilde{\delta}\,y}{\beta\sqrt{v+d-1+y^2}};0,1,v+d\right)\right) * 2G_1\left(\frac{1}{\beta^2}\sqrt{\tilde{\delta}'\tilde{\delta}}\,y\sqrt{\frac{(v+1)}{v+y^2}};0,1,v+1\right)\right\}$

$= H(X;\mu,S,\delta)$

*Example 1.* As a simple illustrative example, explains the relationship between the parameters α , δ and v with Shannon and Renyi entropies in one, two and three dimension spaces. Consider $X \sim MST_d(\mu,S,\delta,v)$ with the following cases:

Case (1) d=1, $\mu = 0.3$, $S = 1.5$, $\delta = 0.3$

Case (2) d=2, $\mu = \begin{pmatrix}3\\2\end{pmatrix}, S = \begin{pmatrix}0.7 & 0.3\\0.3 & 3\end{pmatrix}, \delta = \begin{pmatrix}0.3\\2\end{pmatrix}$

Case(3) d=3, $\mu = \begin{pmatrix}0\\0\\0\end{pmatrix}, S = \begin{pmatrix}1 & 0 & 0\\0 & 1 & 0\\0 & 0 & 1\end{pmatrix}, \delta = \begin{pmatrix}0.3\\2\\0.3\end{pmatrix}$

Table 1. Summarizes the values of Rényi entropy of $X \sim MST_d(\mu,S,\delta,v)$ such that v=3,4,5,6,8,10,12 , $\alpha = 1,2,3,4,5,6,8,10$ and α converges to finite value . It is shows that, there is a relationship between the value of Rényi

entropy and the values of the parameters α, v and d whereas in case 1. for fixed v=2 the Rényi entropy converges to $R_\alpha(X; \mu, S, \delta, v) = 1.2311$ as $\alpha \to \infty$. Also we noted that the upper bound of Rényi entropy is 1.2311 for any value of v.

**Table 1.** Shannon and Rényi entropies of $MST_d(\mu, S, \delta, v)$ are computed for $\alpha = 2,3,4,5,6,8,10$ and α converges to infinite, v=3,4,5,6,8,10 and v=12 in one, two and three dimensions.

| d | v | Shannon entropy $H(x; \mu, S, \delta, v)$ | $R_\alpha(x; \mu, S, \delta, v)$ | | | | | | | |
|---|---|---|---|---|---|---|---|---|---|---|
| | | | α = 2 | α = 3 | α = 4 | α = 5 | α = 6 | α = 8 | α = 10 | α → ∞ |
| 1 | 3 | 1.9590 | 1.6571 | 1.5380 | 1.4788 | 1.4352 | 1.4053 | 1.3638 | 1.3371 | 1.2311 |
| | 4 | 1.8678 | 1.6033 | 1.5010 | 1.4438 | 1.4043 | 1.3749 | 1.3378 | 1.3127 | 1.2101 |
| | 5 | 1.8130 | 1.5750 | 1.4806 | 1.4214 | 1.3839 | 1.3573 | 1.3199 | 1.2958 | 1.1970 |
| | 6 | 1.7767 | 1.5538 | 1.4624 | 1.4089 | 1.3708 | 1.3462 | 1.3086 | 1.2860 | 1.1972 |
| | 8 | 1.7314 | 1.5361 | 1.4459 | 1.3951 | 1.3557 | 1.3310 | 1.2933 | 1.2726 | 1.1970 |
| | 10 | 1.7025 | 1.5140 | 1.4261 | 1.3795 | 1.3451 | 1.3154 | 1.2896 | 1.2624 | 1.1887 |
| | 12 | 1.6871 | 1.4999 | 1.4223 | 1.3688 | 1.3374 | 1.3132 | 1.2780 | 1.2585 | 1.1746 |
| 2 | 3 | 3.5238 | 2.9363 | 2.6728 | 2.5441 | 2.4486 | 2.3939 | 2.3081 | 2.2532 | 2.0756 |
| | 4 | 3.3404 | 2.8648 | 2.6401 | 2.5096 | 2.4300 | 2.3754 | 2.2902 | 2.2403 | 2.0679 |
| | 5 | 3.2363 | 2.8158 | 2.5956 | 2.4832 | 2.4069 | 2.3548 | 2.2781 | 2.2304 | 2.0521 |
| | 6 | 3.1763 | 2.7936 | 2.5964 | 2.4727 | 2.4012 | 2.3468 | 2.2670 | 2.2183 | 2.0625 |
| | 8 | 3.0826 | 2.7565 | 2.5621 | 2.4494 | 2.3849 | 2.3273 | 2.2570 | 2.2109 | 2.0561 |
| | 10 | 3.0367 | 2.7228 | 2.5526 | 2.4386 | 2.3702 | 2.3263 | 2.2508 | 2.2128 | 2.0541 |
| | 12 | 3.0056 | 2.7069 | 2.5472 | 2.4348 | 2.3709 | 2.3164 | 2.2466 | 2.1993 | 2.0539 |
| 3 | 3 | 4.6973 | 3.6934 | 3.2816 | 3.0665 | 2.9312 | 2.8325 | 2.7018 | 2.6177 | 2.3562 |
| | 4 | 4.4806 | 3.6356 | 3.2594 | 3.0701 | 2.9340 | 2.8508 | 2.7266 | 2.6461 | 2.3886 |
| | 5 | 4.3250 | 3.5950 | 3.2465 | 3.0629 | 2.9416 | 2.8526 | 2.7345 | 2.6568 | 2.4839 |
| | 6 | 4.2433 | 3.5636 | 3.2349 | 3.0662 | 2.9399 | 2.8624 | 2.7452 | 2.6711 | 2.4982 |
| | 8 | 4.1010 | 3.5157 | 3.2376 | 3.0577 | 2.9473 | 2.8719 | 2.7564 | 2.6875 | 2.5103 |
| | 10 | 4.0454 | 3.5004 | 3.2102 | 3.0534 | 2.9437 | 2.8705 | 2.7626 | 2.6845 | 2.5271 |
| | 12 | 3.9951 | 3.4695 | 3.2069 | 3.0517 | 2.9414 | 2.8697 | 2.7597 | 2.6840 | 2.5259 |

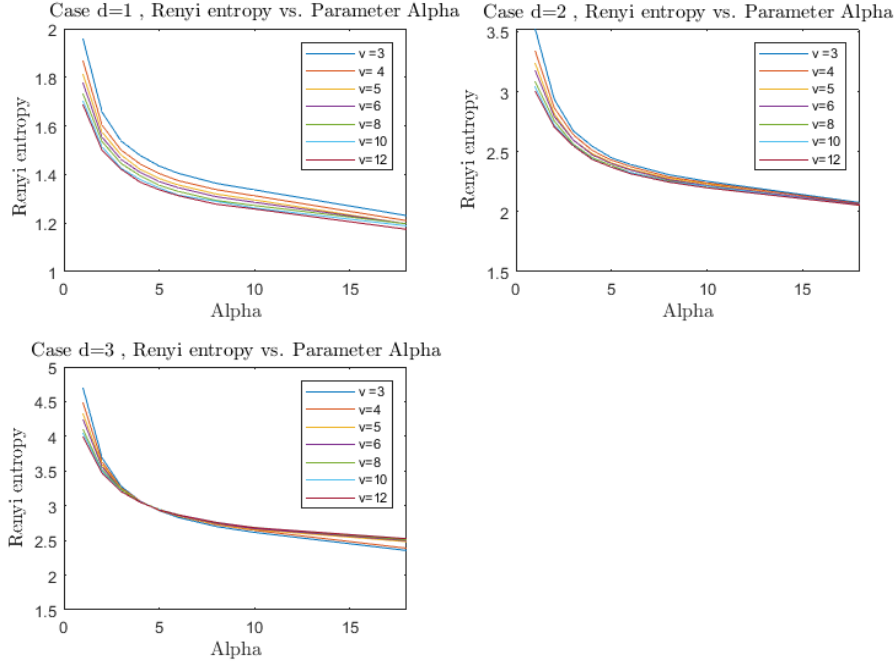

**Figure 1.** The horizontal line represents the values of parameter $\alpha$ and the vertical lines are the Shannon and Rényi entropies of $X \sim MST_d(\mu, S, \delta, v)$ with parameters in example 1.

We observe in Figure 1., that the Rényi entropy of X~$\text{MST}_d(\mu, S, \delta, v)$ converges to finite value for the values of $\alpha$ and v. It can be seen that also, the Shannon and Rényi entropies increases for increases values of d. The dispersion matrices play an important role in determining the value of Rényi entropy where it increases by their increases. Also we note that the effectiveness of $\alpha$ on Rényi entropy is slow whenever the value of $\alpha$ is large and vice versa when it is small, while the effectiveness of v is fickle and it depends on the value of $\alpha$.

## 4. Approximate Shannon and Rényi Entropies for Finite Mixture of Multivariate Skew t- Distributions

**Lemma 14.** Let X~$\text{FMMST}_d(\mu, S, \delta, v, \varepsilon)$. Then the following inequality is accomplished

$$C_{lower} \leq H(X; \mu, S, \delta, v, \varepsilon) \leq C_{upper} \tag{29}$$

where,

$$C_{upper} = \frac{1}{2}\ln(\det(2\pi\exp(1)R)) \tag{30}$$

$$C_{lower} = \sum_{i=1}^{m} \varepsilon_i \left( \ln\left( \frac{\Gamma(\frac{v_i}{2})(v_i\pi)^{\frac{d}{2}}}{\Gamma(\frac{v_i+d}{2})} \right) + \frac{1}{2}\ln(\det(S_i)) + \left(\frac{v_i+d}{2}\right)\{\psi(v_i+d) - \psi(v_i)\} - \mathcal{K}_i \right) \tag{31}$$

$$\mathcal{K}_i = E\left\{ \ln\left( 2G_1\left( \frac{\sqrt{(v_i+d)\tilde{\delta}_i'\tilde{\delta}_i} \, y}{\beta_i\sqrt{v_i+d-1+y^2}}; 0,1, v_i+d \right) \right) * 2G_1\left( \frac{1}{\beta_i^2}\sqrt{\tilde{\delta}_i'\tilde{\delta}_i} \, y\sqrt{\frac{(v_i+1)}{v_i+y^2}}; 0,1, v_i+d \right) \right\} \tag{32}$$

$$R = \sum_{i=1}^{m} \varepsilon_i \frac{v_i}{v_i-2} S_i - \left( \sum_{i=1}^{m} \varepsilon_i \left( \mu_i + \frac{\Gamma(\frac{v_i-1}{2})\sqrt{\frac{v_i}{\pi}}}{\Gamma(\frac{v_i}{2})} \hat{\delta}_i \right) \right) \left( \sum_{i=1}^{m} \varepsilon_i \left( \mu_i + \frac{\Gamma(\frac{v_i-1}{2})\sqrt{\frac{v_i}{\pi}}}{\Gamma(\frac{v_i}{2})} \hat{\delta}_i \right) \right)'$$

(33)

$y \sim T_1(0,1, v_i + d - 1)$

**proof**

Since R is a covariance matrix of X then by using lemma 1., we get the upper bound

$H(X; \mu, S, \delta, v, \varepsilon) \leq \frac{1}{2}\ln(\det(2\pi\exp(1)R))$

Now, to find the lower bound of Rényi entropy, from the expression of mixture density in (4) the Shannon entropy can be written in the following form

$H(X; \mu, S, \delta, v, \varepsilon) = -E\left( \ln\left(\sum_{i=1}^{m} \varepsilon_i \, f(x; \mu_i, S_i, \delta_i, v_i)\right) \right)$

Since $\ln(x)$ is a convex function, then $-\ln(x)$ is a concave. Therefore, by using Jensen's inequality, we obtain

$H(X; \mu, S, \delta, v, \varepsilon) \geq \sum_{i=1}^{m} \varepsilon_i \, H(X; \mu_i, S_i, \delta_i, v_i)$

the Shannon entropy of each component is obtained from proposition 11. This complete the proof.

**Lemma 15.** If X~$\text{FMMST}_d(\mu, S, \delta, v, \varepsilon)$, then for any positive integer $\alpha$ the following inequality is hold

$$R_\alpha(X; \mu, S, \delta, v, \varepsilon) \leq \mathfrak{C}_{Upper} \tag{34}$$

where,

$$\mathfrak{C}_{Upper} = \frac{1}{1-\alpha} \ln\{ \exp((1-\alpha)R_\alpha(X; \mu_m, S_m, \delta_m, v_m)) + \sum_{i=1}^{m-1} \left(\sum_{k=1}^{i} \varepsilon_k\right)^\alpha \left( \exp((1-\alpha)R_\alpha(X; \mu_i, S_i, \delta_i, v_i)) - \exp((1-\alpha)R_\alpha(X; \mu_{i+1}, S_{i+1}, \delta_{i+1}, v_{i+1})) \right) \} \tag{35}$$

**Proof**

Finite mixture density in (4) implies that $\left(f(x; \mu, S, \delta, v, \varepsilon)\right)^\alpha = \left(\sum_{i=1}^{m} \varepsilon_i \, f(x; \mu_i, S_i, \delta_i, v_i)\right)^\alpha$ Using lemma 7. when $p = \alpha$, we get,

$\left(\sum_{i=1}^{m} \varepsilon_i \, f(x; \mu_i, S_i, \delta_i, v_i)\right)^\alpha$
$\geq f(x; \mu_m, S_m, \delta_m, v_n)^\alpha + \sum_{i=1}^{m-1} \left(\sum_{k=1}^{i} \varepsilon_k\right)^\alpha \left( \left(f(x; \mu_i, S_i, \delta_i, v_i)\right)^\alpha - \left(f(x; \mu_{i+1}, S_{i+1}, \delta_{i+1}, v_{i+1})\right)^\alpha \right)$

Taking the integral for both sides of above inequality over $R^d$, we obtained.

$\int_{R^d} (f(x; \mu, S, \delta, v, \varepsilon))^\alpha dx$
$\geq \int_{R^d} f(x; \mu_m, S_m, \delta_m, v_m)^\alpha dx$
$+ \sum_{i=1}^{m-1} (\sum_{k=1}^{i} \varepsilon_k)^\alpha \int_{R^d} [(f(x; \mu_i, S_i, \delta_i, v_i))^\alpha - (f(x; \mu_{i+1}, S_{i+1}, \delta_{i+1}, v_{i+1}))^\alpha] dx$

Again, taking natural logarithm and multiplying by $\frac{1}{1-\alpha}$ for both sides of last inequality, we have
$R_\alpha(X; \mu, S, \delta, v, \varepsilon)$
$\leq \frac{1}{1-\alpha} \ln \{ \int_{R^d} f(x; \mu_m, S_m, \delta_m, v_m)^\alpha dx$
$+ \sum_{i=1}^{m-1} (\sum_{k=1}^{i} \varepsilon_k)^\alpha \int_{R^d} [(f(x; \mu_i, S_i, \delta_i, v_i))^\alpha - (f(x; \mu_{i+1}, S_{i+1}, \delta_{i+1}, v_{i+1}))^\alpha] dx \}$

**Lemma 16.** Let $X \sim FMMST_d(\mu, S, \delta, v, \varepsilon)$. Then for any positive integers $k_1, k_2, \ldots, k_m$ such that $\sum_{i=1}^{m} k_i = \alpha$ the following approximation

$$\frac{1}{\alpha} \ln \left\{ \left( \frac{\alpha!}{k_1! \, k_2! \ldots k_m!} \right) \prod_{i=1}^{m} (\varepsilon_i f(x; \mu_i, S_i, \delta_i, v_i))^{k_i} \right\} \cong -\sum_{i=1}^{m} \gamma_i \ln \left( \frac{\gamma_i}{\varepsilon_i f(x; \mu_i, S_i, \delta_i, v_i)} \right) \tag{36}$$

is accomplished as $\alpha \to \infty$, where $\gamma_i = \frac{k_i}{\alpha}, i = 1,2, \ldots, m$

**Proof**
$\frac{1}{\alpha} \ln \left\{ \left( \frac{\alpha!}{k_1! \, k_2! \ldots k_m!} \right) \prod_{i=1}^{m} (\varepsilon_i f(x; \mu_i, S_i, \delta_i, v_i))^{k_i} \right\}$
$= \frac{1}{\alpha} \ln(\alpha!) - \frac{1}{\alpha} \sum_{i=1}^{m} \ln(k_i!) + \frac{1}{\alpha} \sum_{i=1}^{m} k_i \ln(\varepsilon_i f(x; \mu_i, S_i, \delta_i, v_i))$

By using the approximation of factorial in above equality, we have
$\frac{1}{\alpha} \ln \left\{ \left( \frac{\alpha!}{k_1! \, k_2! \ldots k_m!} \right) \prod_{i=1}^{m} (\varepsilon_i f(x; \mu_i, S_i, \delta_i, v_i))^{k_i} \right\}$
$= \ln(\alpha) - 1 + \frac{1}{2\alpha} \ln(2\pi\alpha) - \frac{1}{\alpha} \sum_{i=1}^{m} k_I \ln(k_i)$
$+ \frac{1}{\alpha} \sum_{i=1}^{m} k_i - \frac{1}{2\alpha} \sum_{i=1}^{m} \ln(2\pi k_i) + \sum_{i=1}^{m} \gamma_i \ln(\varepsilon_i f(x; \mu_i, S_i, \delta_i, v_i))$

But $\gamma_i = \frac{k_i}{\alpha}$, $i = 1,2, \ldots, m$, then $\sum_{i=1}^{m} \gamma_i = 1$.
Consequently,
$\frac{1}{\alpha} \ln \left\{ \left( \frac{\alpha!}{k_1! \, k_2! \ldots k_m!} \right) \prod_{i=1}^{m} (\varepsilon_i f(x; \mu_i, S_i, \delta_i, v_i))^{k_i} \right\} = \frac{1}{2\alpha} \left[ \ln \left( \frac{(2\pi\alpha)^{1-n}}{\prod_{i=1}^{n} \gamma_i} \right) \right] - \sum_{i=1}^{m} \gamma_i \ln \left( \frac{\gamma_i}{\varepsilon_i f(x; \mu_i, S_i, \delta_i, v_i)} \right)$

But $\lim_{\alpha \to \infty} \sum_{i=1}^{m} \gamma_i \ln \left( \frac{\gamma_i}{\varepsilon_i f(x; \theta_i)} \right) = 0$

This implies that the result in this lemma is satisfied.

**Lemma 17.** The approximation
$$R_\alpha(X; \mu, S, \delta, v, \varepsilon) \cong \frac{1}{1-\alpha} \ln \left( \sum_{k_i \in B} (\prod_{i=1}^{m} (\gamma_i)^{-k_i}) \left( \prod_{i=1}^{m} \varepsilon_i^{k_i} \exp \left( (1-\alpha) R_{k_i}(X; \mu_i, S_i, \delta_i, v_i) \right) \right) \right) \tag{37}$$

is satisfied as $\alpha \to \infty$.
where, $\sum_{k_i \in B} \frac{\alpha!}{\prod_{i=1}^{m} k_i!} = m^\alpha$, $B = \{k_i \in N, k_i > 0, \sum_{i=1}^{m} k_i = \alpha, i = 1,2, \ldots, m\}$

**Proof**
The integral of mixture model in (4) with exponent $\alpha$ can be written as
$\int_{R^d} (f(x; \mu, S, \delta, v, \varepsilon))^\alpha dx = \int_{R^d} (\sum_{i=1}^{m} \varepsilon_i f(x; \mu_i, S_i, \delta_i, v_i))^\alpha$
By using multinomial theorem, we obtain
$$\int_{R^d} (f(x; \mu, S, \delta, v, \varepsilon))^\alpha dx = \int_{R^d} \sum_{k_i \in B} \frac{\alpha!}{\prod_{i=1}^{m} k_i!} \prod_{i=1}^{m} (\varepsilon_i f(x; \mu_i, S_i, \delta_i, v_i))^{k_i} dx \tag{38}$$
where, $\sum_{k_i \in B} \frac{\alpha!}{\prod_{i=1}^{m} k_i!} = m^\alpha$, $B = \{k_i \in N, k_i > 0, \sum_{i=1}^{m} k_i = \alpha, i = 1,2, \ldots, m\}$
By replacing right side of equation (38) in equation (36), we get
$\int_{R^d} (f(x; \mu, S, \delta, v, \varepsilon))^\alpha dx \cong \int_{R^d} \sum_{k_i \in B} \exp \left\{ -\alpha \sum_{i=1}^{m} \gamma_i \ln \left( \frac{\gamma_i}{\varepsilon_i f(x; \mu_i, S_i, \delta_i, v_i)} \right) \right\} dx$
The last approximation can be written as
$\int_{R^d} (f(x; \mu, S, \delta, v, \varepsilon))^\alpha = \sum_{k_i \in B} \prod_{i=1}^{m} (\gamma_i)^{-k_i} \int_{R^d} \prod_{i=1}^{m} (\varepsilon_i f(x; \mu_i, S_i, \delta_i, v_i))^{k_i} dx$

Consequently,

$$\int_{\mathbb{R}^d}(f(x;\mu,S,\delta,v,\varepsilon))^\alpha = \sum_{k_i \in B}\left[\prod_{i=1}^m (\gamma_i)^{-k_i}\right]\left[\prod_{i=1}^m \int_{\mathbb{R}^d}(\varepsilon_i f(x;\mu_i,S_i,\delta_i,v_i))^{k_i} dx\right]$$

If we take the natural logarithm and multiplying by $\frac{1}{1-\alpha}$ for both sides of last approximation, then the proof is completed.

**Lemma 18.** Consider $X \sim FMMST_d(\mu,S,\delta,v,\varepsilon)$, then for any positive integer $\alpha$ the lower bound of Rényi entropy appears as follow

$$R_\alpha(X;\mu,S,\delta,v,\varepsilon) \geq \mathfrak{C}_{Lower}$$

where, $\mathfrak{C}_{Lower} = \frac{1}{1-\alpha}\ln\left(\sum_{k_i \in B}\frac{\alpha!}{\prod_{i=1}^m k_i!}\prod_{i=1}^m (\varepsilon_i)^{k_i}\exp\left\{\frac{(1-\alpha)}{\alpha}\sum_{i=1}^m k_i R_\alpha(X;\mu_i,S_i,\delta_i,v_i)\right\}\right)$  (39)

**Proof**

The Rényi entropy of $X \sim FMMST_d(\mu,S,\delta,v,\varepsilon)$ can be written in the following form

$$R_\alpha(X;\mu,S,\delta,v,\varepsilon) = \frac{1}{1-\alpha}\ln\left(\int_{\mathbb{R}^d}\left(\sum_{i=1}^m \varepsilon_i f(x;\mu_i,S_i,\delta_i,v_i)\right)^\alpha dx\right)$$

By using multinomial theorem, we obtain

$$\int_{\mathbb{R}^d}(f(x;\mu,S,\delta,v,\varepsilon))^\alpha dx = \sum_{k_i \in B}\frac{\alpha!}{\prod_{i=1}^m k_i!}\prod_{i=1}^m (\varepsilon_i)^{k_i}\int_{\mathbb{R}^d}\prod_{i=1}^m (f(x;\mu_i,S_i,\delta_i,v_i))^{k_i} dx$$

Applying Hölder's Inequality for m-products, we get

$$\int_{\mathbb{R}^d}(f(x;\mu,S,\delta,v,\varepsilon))^\alpha dx \leq \sum_{k_i \in B}\frac{\alpha!}{\prod_{i=1}^m k_i!}\prod_{i=1}^m (\varepsilon_i)^{k_i}\prod_{i=1}^m \left(\int_{\mathbb{R}^d}(f(x;\mu_i,S_i,\delta_i,v_i))^{p_i k_i} dx\right)^{\frac{1}{p_i}}$$

Such that $p_1, p_2, \ldots, p_m > 0$ and $\sum_{i=1}^m \frac{1}{p_i} = 1$.

the last equation can be written as

$$\int_{\mathbb{R}^d}(f(x;\mu,S,\delta,v,\varepsilon))^\alpha dx \leq \sum_{k_i \in B}\frac{\alpha!}{\prod_{i=1}^m k_i!}\prod_{i=1}^m (\varepsilon_i)^{k_i}\exp\left\{\sum_{i=1}^m\left(\frac{(1-p_i k_i)}{p_i}R_{p_i k_i}(X;\mu_i,S_i,\delta_i,v_i)\right)\right\}$$

we choose $p_i = \frac{\alpha}{k_i}$, $i = 1,2,\ldots,m$ such that $\sum_{i=1}^m \frac{1}{p_i} = \sum_{i=1}^m \frac{k_i}{\alpha} = 1$ and $1 \leq \frac{\alpha}{k_i} \leq \alpha$, to obtain

$$\int_{\mathbb{R}^d}(f(x;\mu,S,\delta,v,\varepsilon))^\alpha dx \leq \sum_{k_i \in B}\frac{\alpha!}{\prod_{i=1}^m k_i!}\prod_{i=1}^m (\varepsilon_i)^{k_i}\exp\left\{\frac{(1-\alpha)}{\alpha}\sum_{i=1}^m k_i R_\alpha(X;\mu_i,S_i,\delta_i,v_i)\right\}$$

We complete the proof by taking the natural logarithm and multiplying by $\frac{1}{1-\alpha}$ for both sides of last inequality.

**Theorem 19.** Let $X \sim FMMST_d(\mu,S,\delta,v,\varepsilon)$. Then the approximate Renyi entropy of X appears as follow

$$R_\alpha(X;\mu,S,\delta,v,\varepsilon) = \frac{1}{2(1-\alpha)}\left\{\ln\left(\sum_{k_i \in B}\frac{\alpha!}{\prod_{i=1}^m k_i!}\prod_{i=1}^m (\varepsilon_i)^{k_i}\exp\left\{\frac{(1-\alpha)}{\alpha}\sum_{i=1}^m k_i R_\alpha(X;\mu_i,S_i,\delta_i,v_i)\right\}\right)\right.$$
$$+ (1-\alpha)R_\alpha(X;\mu_m,S_m,\delta_m,v_m) + \ln\left(\sum_{i=1}^{m-1}\left(\sum_{k=1}^i \varepsilon_k\right)^\alpha\right.$$
$$\left.\left.\cdot\left(\exp\left((1-\alpha)R_\alpha(X;\mu_i,S_i,\delta_i,v_i)\right) - \exp\left((1-\alpha)R_\alpha(X;\mu_{i+1},S_{i+1},\delta_{i+1},v_{i+1})\right)\right)\right)\right\}$$

**Proof**

The proof is directed from lemmas 15. and 18., by taking the mean of upper and lower bounds.

*Example 2.* As a simple illustrative example, explains the relationship between the parameters $\alpha$, $\delta$ and $v$ with Shannon and Renyi entropies in one, two and three dimension spaces. Consider $FMMST_d(\mu,S,\delta,v,\varepsilon)$ with the following cases:

**Case (1):** d=1
m=2, $\varepsilon = (0.2, 0.8)$,  $\mu = (0.3, 4)$,  S = (1.5, 5)  $\delta = (0.3, 4)$  and v=(3,3)
m=3, $\varepsilon = (0.2, 0.3, 0.5)$,  $\mu = (0.3, 4, 0.6)$,  S = (1.5, 5, 3), $\delta = (0.3, 4, 2.2)$ and v=(3,3,4)
m=4, $\varepsilon = (0.1, 0.2, 0.2, 0.5)$, $\mu = (0.3, 4, 0.6, 3)$, S = (1.5, 5, 3, 2), $\delta = (0.3, 4, 2.2, 1)$, v=(3,3,4,4)
m=5, $\varepsilon = (0.2, 0.2, 0.2, 0.2, 0.2)$, $\mu = (0.3, 4, 0.6, 3, 2)$, S = (1.5, 5, 3, 2, 5),
    $\delta = (0.3, 4, 2.2, 1, 2.1)$ and v=(3,3,4,4,5)

**Case (2):** d=2
m=2, $\varepsilon = (0.2, 0.8)$, $\mu = \left(\binom{3}{2}, \binom{1}{5}\right)$, $S = \left(\begin{pmatrix}0.7 & 0.3\\0.3 & 3\end{pmatrix}, \begin{pmatrix}0.12 & 0.13\\0.13 & 3\end{pmatrix}\right)$, $\delta = \left(\binom{0.16}{0.59}, \binom{2.3}{3.1}\right)$ and v=(3,3)

m=3 , ε = (0.2,0.3,0.5) , μ = $\left(\binom{3}{2},\binom{1}{5},\binom{3}{1}\right)$ , S = $\left(\begin{pmatrix}0.7 & 0.3\\0.3 & 3\end{pmatrix},\begin{pmatrix}0.12 & 0.13\\0.13 & 3\end{pmatrix},\begin{pmatrix}0.18 & 0.6\\0.6 & 4\end{pmatrix}\right)$ ,

δ = $\left(\binom{0.16}{0.59},\binom{2.3}{3.1},\binom{2.6}{1}\right)$ and v=(3,3,4)

m=4 , ε = (0.1,0.2,0.2,0.5), μ = $\left(\binom{3}{2},\binom{1}{5},\binom{3}{1},\binom{1}{1}\right)$,

S = $\left(\begin{pmatrix}0.7 & 0.3\\0.3 & 3\end{pmatrix},\begin{pmatrix}0.12 & 0.13\\0.13 & 3\end{pmatrix},\begin{pmatrix}0.18 & 0.6\\0.6 & 4\end{pmatrix},\begin{pmatrix}1 & 0\\0 & 1\end{pmatrix}\right)$,

δ = $\left(\binom{0.16}{0.59},\binom{2.3}{3.1},\binom{2.6}{1},\binom{0.6}{1}\right)$ and v=(3,3,4,4)

m=5 , ε = (0.2,0.2,0.2,0.2,0.2), μ = $\left(\binom{3}{2},\binom{1}{5},\binom{3}{1},\binom{1}{1},\binom{1}{0.3}\right)$,

S = $\left(\begin{pmatrix}0.7 & 0.3\\0.3 & 3\end{pmatrix},\begin{pmatrix}0.12 & 0.13\\0.13 & 3\end{pmatrix},\begin{pmatrix}0.18 & 0.6\\0.6 & 4\end{pmatrix},\begin{pmatrix}1 & 0\\0 & 1\end{pmatrix},\begin{pmatrix}1 & 0\\0 & 1\end{pmatrix}\right)$,

δ = $\left(\binom{0.16}{0.59},\binom{2.3}{3.1},\binom{2.6}{1},\binom{0.6}{1},\binom{1}{1}\right)$ and v=(3,3,4,4,5)

**Case (3) :** d=3

m=2, ε = (0.2,0.8), μ = $\left(\begin{pmatrix}3\\2\\3\end{pmatrix},\begin{pmatrix}1\\5\\1\end{pmatrix}\right)$, S = $\left(\begin{pmatrix}0.7 & 0.3 & 0.5\\0.3 & 3 & 0.3\\0.5 & 0.3 & 1\end{pmatrix},\begin{pmatrix}5 & 0.3 & 2\\0.3 & 5 & 1\\2 & 1 & 3\end{pmatrix}\right)$

, δ = $\left(\begin{pmatrix}0.16\\0.59\\0.1\end{pmatrix},\begin{pmatrix}2.3\\3.1\\1.5\end{pmatrix}\right)$ and v=(3,3)

m=3 , ε = (0.2,0.3,0.5), μ = $\left(\begin{pmatrix}3\\2\\3\end{pmatrix},\begin{pmatrix}1\\5\\1\end{pmatrix},\begin{pmatrix}2\\3\\1\end{pmatrix}\right)$,

S = $\left(\begin{pmatrix}0.7 & 0.3 & 0.5\\0.3 & 3 & 0.3\\0.5 & 0.3 & 1\end{pmatrix},\begin{pmatrix}5 & 0.3 & 2\\0.3 & 5 & 1\\2 & 1 & 3\end{pmatrix},\begin{pmatrix}1 & 0 & 0\\0 & 1 & 0\\0 & 0 & 1\end{pmatrix}\right)$, δ = $\left(\begin{pmatrix}0.16\\0.59\\0.1\end{pmatrix},\begin{pmatrix}2.3\\3.1\\1.5\end{pmatrix},\begin{pmatrix}2\\1\\2\end{pmatrix}\right)$ and v=(3,3,4)

m=4 , ε = (0.1,0.2,0.2,0.5), μ = $\left(\begin{pmatrix}3\\2\\3\end{pmatrix},\begin{pmatrix}1\\5\\1\end{pmatrix},\begin{pmatrix}2\\3\\1\end{pmatrix},\begin{pmatrix}1\\1\\1\end{pmatrix}\right)$,

S = $\left(\begin{pmatrix}0.7 & 0.3 & 0.5\\0.3 & 3 & 0.3\\0.5 & 0.3 & 1\end{pmatrix},\begin{pmatrix}5 & 0.3 & 2\\0.3 & 5 & 1\\2 & 1 & 3\end{pmatrix},\begin{pmatrix}1 & 0 & 0\\0 & 1 & 0\\0 & 0 & 1\end{pmatrix},\begin{pmatrix}1 & 0 & 0\\0 & 1 & 0\\0 & 0 & 1\end{pmatrix}\right)$ , δ = $\left(\begin{pmatrix}0.16\\0.59\\0.1\end{pmatrix},\begin{pmatrix}2.3\\3.1\\1.5\end{pmatrix},\begin{pmatrix}2\\1\\2\end{pmatrix},\begin{pmatrix}1\\2\\1\end{pmatrix}\right)$ and v=(3,3,4,4)

m=5 , ε = (0.2,0.2,0.2,0.2,0.2), μ = $\left(\begin{pmatrix}3\\2\\3\end{pmatrix},\begin{pmatrix}1\\5\\1\end{pmatrix},\begin{pmatrix}2\\3\\1\end{pmatrix},\begin{pmatrix}1\\1\\1\end{pmatrix},\begin{pmatrix}0.1\\1\\0.1\end{pmatrix}\right)$,

S = $\left(\begin{pmatrix}0.7 & 0.3 & 0.5\\0.3 & 3 & 0.3\\0.5 & 0.3 & 1\end{pmatrix},\begin{pmatrix}5 & 0.3 & 2\\0.3 & 5 & 1\\2 & 1 & 3\end{pmatrix},\begin{pmatrix}1 & 0 & 0\\0 & 1 & 0\\0 & 0 & 1\end{pmatrix},\begin{pmatrix}1 & 0 & 0\\0 & 1 & 0\\0 & 0 & 1\end{pmatrix},\begin{pmatrix}1 & 0 & 0\\0 & 1 & 0\\0 & 0 & 1\end{pmatrix}\right)$,

δ = $\left(\begin{pmatrix}0.16\\0.59\\0.1\end{pmatrix},\begin{pmatrix}2.3\\3.1\\1.5\end{pmatrix},\begin{pmatrix}2\\1\\2\end{pmatrix},\begin{pmatrix}1\\2\\1\end{pmatrix},\begin{pmatrix}2\\1\\2\end{pmatrix}\right)$ and v=(3,3,4,4,5)

**Table 2.** Shannon entropy of X∼FMMST$_d$(μ, S, δ, ν, ε) is computed for m = 2,3,4 and 5 in one, two and three dimensions.

| Case | d | m | Approximate Shannon entropy | | | error |
|---|---|---|---|---|---|---|
| | | | $C_{lower}$ | $C_{upper}$ | H | Less than |
| **1** | 1 | 2 | 2.0984 | 2.5555 | 2.3283 | 0.2262 |
| | 1 | 3 | 1.9818 | 2.3523 | 2.1671 | 0.1852 |
| | 1 | 4 | 1.9398 | 2.2569 | 2.0983 | 0.1585 |
| | 1 | 5 | 1.9471 | 2.2569 | 2.1020 | 0.1549 |
| **2** | 2 | 2 | 2.6149 | 3.6332 | 3.1240 | 0.5091 |
| | 2 | 3 | 2.5553 | 3.5219 | 3.0386 | 0.4833 |
| | 2 | 4 | 2.7972 | 3.6936 | 3.2454 | 0.4482 |
| | 2 | 5 | 2.8443 | 3.7310 | 3.2876 | 0.4433 |
| **3** | 3 | 2 | 6.3749 | 7.4959 | 6.9354 | 0.5605 |
| | 3 | 3 | 5.2607 | 6.6252 | 5.9429 | 0.6822 |
| | 3 | 4 | 4.9670 | 6.2674 | 5.6172 | 0.6502 |
| | 3 | 5 | 5.0167 | 6.3223 | 5.6695 | 0.6528 |

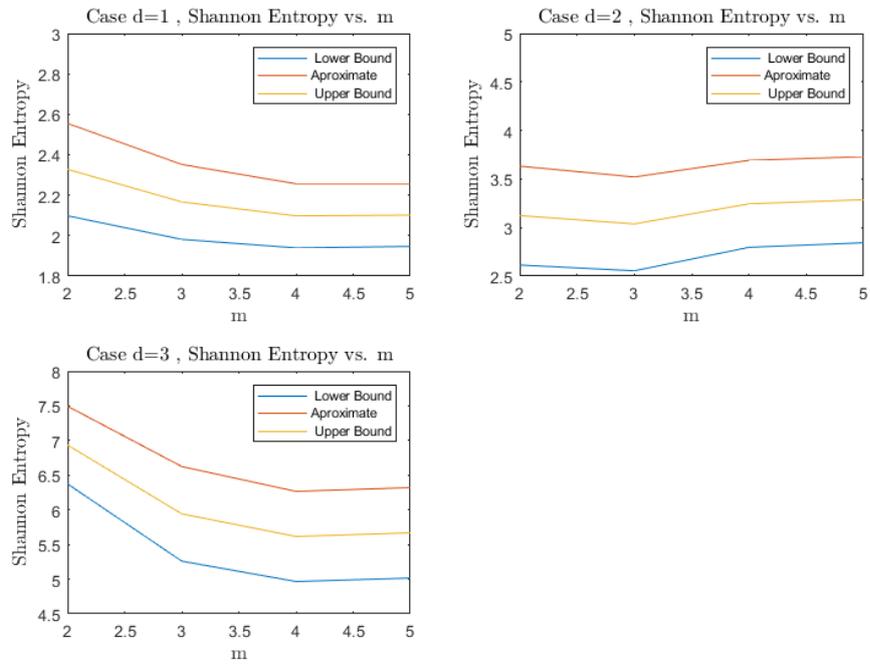

**Figure 2.** The horizontal line represents the values of parameter $m$ and the vertical lines is a Shannon entropy of X∼FMMST$_d$(μ, S, δ, ν, ε) with parameters in example 2.

**Table 3.** Rényi entropy of $X \sim \text{FMST}_d(\mu, S, \delta, v)$ is computed for $\alpha = 2,3,4,5,6,8,10.15,20,30$ and $m=2,3,4$ in one, two and three dimensions.

| Case | d | m | α | Renyi entropy | | Error | |
|---|---|---|---|---|---|---|---|
| | | | | $C_\alpha$ Lower | $C_\alpha$ Upper | Approximate | Less than |
| 1 | 1 | 2 | 2 | 1.8936 | 1.9287 | 1.9112 | 0.0176 |
| | 1 | 2 | 3 | 1.7668 | 1.8306 | 1.7987 | 0.0319 |
| | 1 | 2 | 4 | 1.7092 | 1.7727 | 1.7410 | 0.0317 |
| | 1 | 2 | 5 | 1.6620 | 1.7263 | 1.6942 | 0.0321 |
| | 1 | 2 | 10 | 1.5659 | 1.6301 | 1.5980 | 0.0321 |
| | 1 | 2 | 15 | 1.5269 | 1.5943 | 1.5606 | 0.0337 |
| | 1 | 2 | 20 | 1.5032 | 1.5687 | 1.5359 | 0.0327 |
| | 1 | 2 | 30 | 1.4784 | 1.5438 | 1.5111 | 0.0327 |
| | 1 | 3 | 2 | 1.7719 | 1.7775 | 1.7747 | 0.0028 |
| | 1 | 3 | 3 | 1.6339 | 1.6668 | 1.6504 | 0.0164 |
| | 1 | 3 | 4 | 1.5739 | 1.6053 | 1.5896 | 0.0157 |
| | 1 | 3 | 5 | 1.5285 | 1.5637 | 1.5461 | 0.0176 |
| | 1 | 3 | 10 | 1.4372 | 1.4698 | 1.4535 | 0.0163 |
| | 1 | 3 | 15 | 1.4001 | 1.4309 | 1.4155 | 0.0154 |
| | 1 | 3 | 20 | 1.3815 | 1.4115 | 1.3965 | 0.0150 |
| | 1 | 3 | 30 | 1.3768 | 1.4066 | 1.3918 | 0.0148 |
| | 1 | 4 | 2 | 1.7063 | 1.7259 | 1.7161 | 0.0098 |
| | 1 | 4 | 3 | 1.5721 | 1.6173 | 1.5947 | 0.0226 |
| | 1 | 4 | 4 | 1.5012 | 1.5569 | 1.5291 | 0.0278 |
| | 1 | 4 | 5 | 1.4687 | 1.5166 | 1.4927 | 0.0239 |
| | 1 | 4 | 10 | 1.3676 | 1.4231 | 1.3954 | 0.0278 |
| | 1 | 4 | 15 | 1.2968 | 1.3577 | 1.3273 | 0.0304 |
| | 1 | 4 | 20 | 1.2614 | 1.3250 | 1.2933 | 0.0317 |
| | 1 | 4 | 30 | 1.2526 | 1.3168 | 1.2848 | 0.0320 |

| Case | d | m | α | Renyi entropy | | Error | |
|---|---|---|---|---|---|---|---|
| | | | | $C_\alpha$ Lower | $C_\alpha$ Upper | Approximate | Less than |
| 2 | 2 | 2 | 2 | 1.7185 | 1.9222 | 1.8204 | 0.1018 |
| | 2 | 2 | 3 | 1.4530 | 1.6230 | 1.5380 | 0.0850 |
| | 2 | 2 | 4 | 1.2991 | 1.4887 | 1.3939 | 0.0948 |
| | 2 | 2 | 5 | 1.2147 | 1.4105 | 1.3126 | 0.0979 |
| | 2 | 2 | 10 | 1.0119 | 1.1883 | 1.1001 | 0.0882 |
| | 2 | 2 | 15 | 0.9275 | 1.1037 | 1.0156 | 0.0881 |
| | 2 | 2 | 20 | 0.8836 | 1.0517 | 0.9676 | 0.0840 |
| | 2 | 2 | 30 | 0.8327 | 1.0034 | 0.9180 | 0.0854 |
| | 2 | 3 | 2 | 1.6939 | 1.8950 | 1.7945 | 0.1005 |
| | 2 | 3 | 3 | 1.4300 | 1.6355 | 1.5327 | 0.1028 |
| | 2 | 3 | 4 | 1.3101 | 1.4966 | 1.4033 | 0.0932 |
| | 2 | 3 | 5 | 1.2127 | 1.3976 | 1.3051 | 0.0925 |
| | 2 | 3 | 10 | 1.0111 | 1.1857 | 1.0984 | 0.0873 |
| | 2 | 3 | 15 | 0.9282 | 1.1005 | 1.0143 | 0.0862 |
| | 2 | 3 | 20 | 0.8849 | 1.0536 | 0.9692 | 0.0843 |
| | 2 | 3 | 30 | 0.8337 | 1.0019 | 0.9178 | 0.0841 |
| | 2 | 4 | 2 | 2.1855 | 2.2517 | 2.2186 | 0.0331 |
|   | 2 | 4 | 3  | 1.9369 | 2.0710 | 2.0040 | 0.0671 |
|   | 2 | 4 | 4  | 1.7970 | 1.9845 | 1.8908 | 0.0937 |
|   | 2 | 4 | 5  | 1.7047 | 1.9253 | 1.8150 | 0.1103 |
|   | 2 | 4 | 10 | 1.6216 | 1.8720 | 1.7468 | 0.1252 |
|   | 2 | 4 | 15 | 1.5634 | 1.8347 | 1.6991 | 0.1356 |
|   | 2 | 4 | 20 | 1.5343 | 1.8160 | 1.6753 | 0.1407 |
|   | 2 | 4 | 30 | 1.5270 | 1.8113 | 1.6694 | 0.1419 |
| 3 | 3 | 2 | 2  | 5.2280 | 5.5015 | 5.3647 | 0.1367 |
|   | 3 | 2 | 3  | 4.7988 | 5.1659 | 4.9823 | 0.1835 |
|   | 3 | 2 | 4  | 4.5649 | 4.9848 | 4.7749 | 0.2100 |
|   | 3 | 2 | 5  | 4.4164 | 4.8589 | 4.6377 | 0.2212 |
|   | 3 | 2 | 10 | 4.0815 | 4.5706 | 4.3260 | 0.2445 |
|   | 3 | 2 | 15 | 3.9489 | 4.4481 | 4.1985 | 0.2496 |
|   | 3 | 2 | 20 | 3.8771 | 4.3877 | 4.1324 | 0.2553 |
|   | 3 | 2 | 30 | 3.8592 | 4.3726 | 4.1159 | 0.2567 |
|   | 3 | 3 | 2  | 3.8026 | 4.1182 | 3.9604 | 0.1578 |
|   | 3 | 3 | 3  | 3.2605 | 3.7012 | 3.4809 | 0.2203 |
|   | 3 | 3 | 4  | 3.0183 | 3.4761 | 3.2472 | 0.2289 |
|   | 3 | 3 | 5  | 2.8843 | 3.3352 | 3.1098 | 0.2254 |
|   | 3 | 3 | 10 | 2.5761 | 3.0129 | 2.7945 | 0.2184 |
|   | 3 | 3 | 15 | 2.4634 | 2.8866 | 2.6750 | 0.2116 |
|   | 3 | 3 | 20 | 2.4070 | 2.8235 | 2.6152 | 0.2082 |
|   | 3 | 3 | 30 | 2.3929 | 2.8077 | 2.6003 | 0.2074 |
|   | 3 | 4 | 2  | 3.6630 | 3.9208 | 3.7919 | 0.1289 |
|   | 3 | 4 | 3  | 3.2517 | 3.5226 | 3.3872 | 0.1355 |
|   | 3 | 4 | 4  | 3.0487 | 3.3094 | 3.1790 | 0.1304 |
|   | 3 | 4 | 5  | 2.9208 | 3.1719 | 3.0464 | 0.1255 |
|   | 3 | 4 | 10 | 2.6246 | 2.8607 | 2.7426 | 0.1180 |
|   | 3 | 4 | 15 | 2.4173 | 2.6429 | 2.5299 | 0.1130 |
|   | 3 | 4 | 20 | 2.3136 | 2.5340 | 2.4236 | 0.1104 |
|   | 3 | 4 | 30 | 2.2877 | 2.5068 | 2.3970 | 0.1098 |

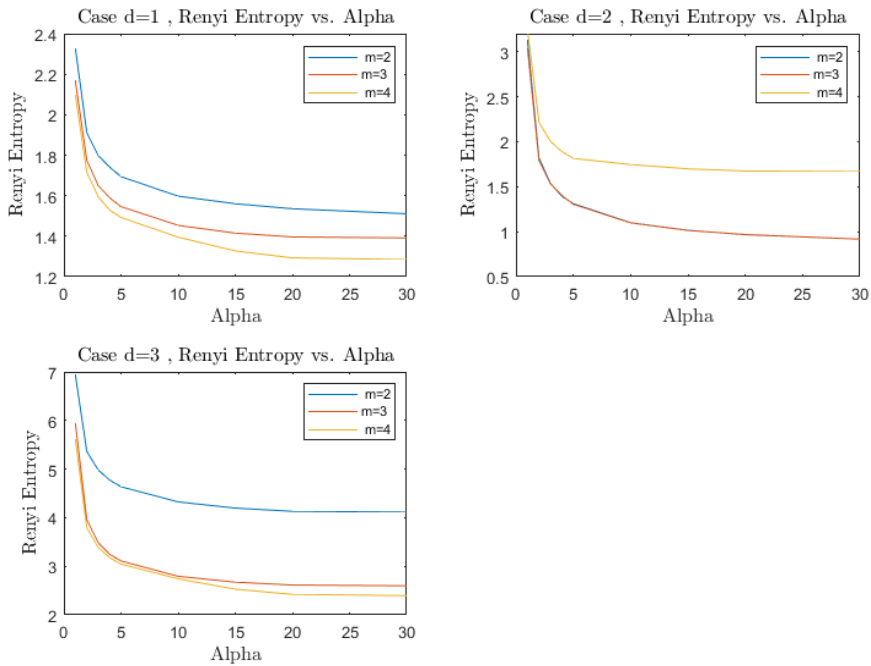

**Figure 3.** The horizontal line represents the values of parameter $\alpha$ and the vertical lines is a Renyi entropy of x~FMMST$_d$(μ, S, δ, v, ε) with parameters in example 2.

## 5. Conclusion

We have derived the lower and upper bounds on the Shannon and Rényi entropies of X~MST$_d$(μ, S, δ, v) and we extended these tools to the class of mixtures models. Using the mean of these bounds, the approximate value of entropy can be calculated. the entropy both types (Shannon and Rényi) converges to finite value of X~MST$_d$(μ, S, δ, v) and its mixtures model for any values of α , v , m and d.

## References

[1] Alan G. and Frank B. (2009). *Computation of Multivariate Normal and t Probabilities*. Springer-Verlag Berlin Heidelberg.
[2] Pav S. (2015). *Moments of the Log non-central Chi-square distribution*. arXiv :1503.06266v1,21 Mar.2015.
[3] Lin P. (1972). *Some Characterization of the Multivariate t Distribution*. Journal of Multivariate Analysis.
[4] Pyne S., Hu X., Wang K., Rossin E., Lin T., Maier L., Baecher-Allan C., McLachlan G., Tamayo P., Hafler D., and Mesirov J. (2009). *Automated High-dimensional Flow Cytometry Data Analysis*. Springer.
[5] Azzalini A. and Regoli G. (2012). *Some Properties of Skew-symmetric Distributions*. Annals of Institute of Statistical Mathematics, Vol. 64,Issue 4,pp 857-879.
[6] Casarin R. (2013). Monte Carlo Methods using Matlab. University of Brescia.
[7] Azzalini A. and Capitanio A. (2003). *Distribution Generated by Perturbation of Symmetry with Emphasis on a Multivariate Skew t Distribution*. Journal of the Royal Statistical Society: series B, 65 (2), 367-389.
[8] José C. (2010). *Principe Information Theoretic Learning Renyi's Entropy and Kernal Perspectives*. Springer Science +Business Media, LLC.
[9] Lee S. and McLachlan G. (2014). *Finite Mixtures of Multivariate Skew t-distributions*. Statistics ana Computing, Dol:10.1007/s11222-012-9362-4.
[10] Contreras-Reyes J. and Cortés D. (2016). *Bounds on Rényi and Shannon Entropies for Finite Mixtures of Multivariate Skew-Normal Distributions: Applications to Swordfish (Xiphias gladius Linnaeus )*. Entropy 11, 382; doi:10.3390/e18110382.
[11] Arellano-Valle R., Contreras-Reyes J. and Genton (2012). *Shannon Entropy and Mutual information for Multivariate Skew Elliptical Distributions*. Scandinavian Journal of Statistics Theory and Applications dio:10.1111/j.1467-9469.2011.00774.
[12] Cover T. M. Thomas (2006). *Elements of Information Theory*. Wily and Son, 2$^{nd}$ New York. NY, USA.
[13] Bennett G.. (1986). *Lower Bounds for Matrices*. Elsevier Science Publishing Co.

[14] Abramowitz M. and Stegun A. (1972). *Psi (Diagamma) Function*. Handbook of Mathematical Functions with Formulas, Graphs and Mathematical Tables (10$^{th}$ ed.) New York.

[15] Javier E. and Contreras-Reyes J. (2016). *Rényi Entropy and Complexity Measure for Skew-Gaussian Distributions and Related Families*. arXiv:1406.0111v2[physics. Data-an] 8 May.

[16] Lin T., Lee J. and Wan H. (2007). *Robust Mixture Modeling Using the Skew t-Distribution*. Statistics and Computing Dol 10.1007/s11222-006-9005-8.

[17] Pyne S., Hu X, Wang K., Rossin E., Lin T., Maier L., Baecher-Allan C., McLachlan G., Tamayo P., Halfer D., De J., Mesirow J. (2009). *Automated High-dimensional Flow Cytometric data Analysis*. Proceedings of the National Academy of Sciences USA 106:8519-8524.

[18] C. E. Shanon (1948). *A mathematical Theory of Communication*. Bell systems technology, j. vol. 27, pp. 379-423.

[19] Rényi A. (1961). *On Measures of Information and entropy. Proceeding of the Fourth Berkeley Symposium on Mathematics, Statistics and Probability pp. 547-561*.

[20] Wood R., Blythe R. and Evans M. (2017). *Rényi Entropy of the Totally Asymmetric Exclusion Process*. arXiv:1708.00303.